\documentclass{amsart}
 \usepackage{amssymb}
\usepackage{amscd}
\usepackage{tikz-cd}
\usepackage{mathtools}
\input xy
\xyoption{all}
\newtheorem{theorem}{Theorem}[section]
\newtheorem{lemma}[theorem]{Lemma}

\theoremstyle{definition}
\newtheorem{definition}[theorem]{Definition}
\newtheorem{example}[theorem]{Example}

\theoremstyle{remark}
\newtheorem{remark}[theorem]{Remark}

\theoremstyle{notation}

\theoremstyle{proposition}
\newtheorem{proposition}[theorem]{Proposition}

\theoremstyle{corollary}
\newtheorem{corollary}[theorem]{Corollary}

\numberwithin{equation}{section}

\begin{document}

\title{On k\"ahler differentials of divided powers algebras}

% Information for first author
%\author{ }
%\address { }
%Current address
%\curraddr{}
%\email{}
     %\thanks will become a 1st page footnote.
%\thanks{}

%Information for second author
\author{Ioannis Dokas}
\address{National Kapodistrian University of Athens}
\email{iodokas@math.uoa.gr}

%    General info
\subjclass[2010]{13D03, 18G10, 17B50, 14A15, 14G17}

%\date{ }

%\dedicatory{}

\keywords{Divided powers algebras, Quillen-Barr-Beck cohomology, K\"ahler differentials}

\begin{abstract} 
The Quillen-Barr-Beck cohomology of augmented algebras with divided powers is defined as the derived functor of Beck derivations. The main theorem of this paper states that  the K\"ahler differentials of an augmented algebra with divided powers in prime characteristic represents Beck derivations. We give a geometrical interpretation of this statement for the sheaf of relative differentials. As an application in modular Lie theory we prove that any special derivation of a divided powers algebra is a Beck derivation and we apply the theorem  to Witt algebras.
\end{abstract}

\maketitle

%% The correct journal style for \specialsection is all uppercase; a known

%% in amsart.cls prevents this, so input must be uppercase until it is

%\specialsection*{This is a Special Section Head}
%\specialsection*{THIS IS A SPECIAL SECTION HEAD}

%%%%%%%%%%%%%%%%%%%%%%%%%%%%%%%%%%%%%%%%%%%%%%%%%%%%%%%%%%%%%%%%%%%%%%%%
%\footnote{Here is an example of a footnote. Notice that this footnote
%text is running on so that it can stand as an example of how a footnote
%with separate paragraphs should be written.
%\par
%And here is the beginning of the second paragraph.}%
%%%%%%%%%%%%%%%%%%%%%%%%%%%%%%%%%%%%%%%%%%%%%%%%%%%%%%%%%%%%%%%%%%%%%%%%

\section*{\textbf{Introduction}}

The Andr\'e-Quillen cohomology of a commutative algebra $R$ with values in a $R$-module $M$ is the Barr-Beck's cotriple cohomology \cite{BB} of $R$ with coefficients in  the derivation functor $Der(-,M)$. Specifically,  in the category of commutative algebras the group of Beck derivations coincides with the usual group of derivations and Beck $R$-modules coincides with $R$-modules. It follows that the $R$-module of K\"ahler differentials represents Beck derivations in the category of commutative algebras. 

A question which arises is what happens in prime characteristic if the commutative algebra is equipped with divided powers. Our motivation comes from the theory of modular Lie algebras. To be more precise, homomorphisms from modular Lie algebras to the module of special derivations of divided powers algebras, are constructed in order to classify simple modular Lie algebras. An essential role for these constructions is played by the divided powers algebra $O(n)$ and the Witt algebra $W_{n}$ of special derivations of $O(n)$. 

To summarize the content of the paper, in Section $1$  we recall some preliminaries from the theory of divided powers algebras. In particular, let $k$ be a field of prime characteristic and  let $A=A_{0}\oplus A_{+}$ be pregraded $k$-algebra (see \cite{Ro2}) with divided powers on the augmentation ideal $A_{+}$. Then we point out the fact that by theorems of H. Cartan \cite{Car} and J.P. Soublin \cite{Sou}  the system divided powers operations on $A_{+}$ are determined by the $p$-th power.

In Section $2$, we recall from \cite{Do}, how is defined Quillen-Barr-Beck cohomology \cite{Q,BB, B} for the category $pdCom_{k}$ of augmented algebras $A=k\oplus A_{+}$ with divided powers on the augmentation ideal $A_{+}$. Precisely, for any  $B\in pdCom_{k}$ divided powers algebra over $A\in pdCom_{k}$ and $M$ a Beck $A$-module is determined in \cite{Do} the group of Beck derivations $Der_{p}(B,M)$ and is proved that the category of Beck $A$-modules is equivalent to the category of modules over a ring $V(A)$. Then the Quillen-Barr-Beck's cohomology $H^{*}(A,M)$ of a divided powers algebra $A\in pdCom_{k}$ with values in a Beck $A$-module $M$ is the cotriple cohomology of $A$ with coefficients in $Der_{p}(-,M)$.

Next we prove the principal result of the paper. Let  $A \in pdCom_{k}$ be a divided powers algebra and let $\mu: A\otimes_{k}A\rightarrow A$ be the multiplication map with $I=ker\,\mu$. We prove that the K\"ahler differentials $I/I^{2}$  has  a richer structure than just an $A$-module. Indeed, the K\"ahler differentials  $I/I^{2}$ is a Beck $A$-module or equivalently a $V(A)$-module. By Theorem $2.9$ we prove that

\begin{equation}
 Hom_{V(A)}(I/I^{2},M)\simeq Der_{p}(A,M).
 \end{equation}
 
Hence the K\"ahler differentials of an augmented algebra $A$ with divided powers on the augmentation ideal in prime characteristic represents Beck derivations. By Proposition $2.3$ we prove that a special derivation of a divided powers algebra in prime characteristic is a Beck derivation. It follows  that the Witt algebra $W_{n}$ is the group of  Beck derivations of $O(n)$.  By Corollary $2.10$ we obtain the bijection
$$H^{0}(O(n),O(n)) \simeq W_{n}.$$
By Proposition $2.11$ we prove that  if $A=k\oplus A_{+}$ is a divided powers algebra and $J$ a is $p$-ideal of $A$ with $J\subset A_{+}$ then the second fundamental exact sequence is an exact sequence of $A/J$-Beck modules.

In the last Section $3$ we give a geometrical interpretation of  the isomorphism $(0.1)$  for the sheaf of relative differentials. Let $A=k\oplus A_{+}$ be a divided powers algebra and let $X=Spec\,A$ be the spectrum of $A$. We prove that the sheaf of differentials $\Delta^{*}\big(\widetilde{I/I^{2}}\big)$ has a richer structure than a $\mathcal{O}_{X}$-module. We call this richer structure $p$-$\mathcal{O}_{X}$-module and we prove by  Theorem $3.13$ that

$$Hom_{p-\mathcal{O}_{X}}( \widetilde{I/I^{2}},\mathcal{M})\simeq Der_{p}(\mathcal{O}_{X},\mathcal{M}),$$ where $Der_{p}(\mathcal{O}_{X},\mathcal{M})$ denotes the set of morphisms of sheaves of $k$-vector spaces $D:\mathcal{O}_{X}\rightarrow \mathcal{M}$ such that $D_{U}$ is a Beck derivation, for all $U\in Open(X)$.

\section{Preliminaries}

H. Cartan in  \cite{Car} proved that  the homotopy of a simplicial commutative algebra is a divided powers algebra.  The theory of divided powers algebras is developed by N. Roby in \cite{Ro1,Ro2,Ro3}.  In the context of crystalline cohomology \cite{Ber}, P. Berthelot  defines the more general notion of ring with divided powers. Let us recall the definition.

\begin{definition}
Let $(A,I)$ be a commutative ring together with an ideal $I\subset A$. A system of divided powers on I,  is a collection of maps  $\gamma_{i}: I\rightarrow A$, for $i\geq 0$, such that the following relations hold

\begin{equation}
\gamma_{0}(x)=1, \; \gamma_{1}(x)=x, \gamma_{i}(x)\in I,\; i\geq 1,
\end{equation}

\begin{align}
\gamma_{i}(x+y) &=\sum_{k=0}^{k=i}\gamma_{k}(x)\gamma_{i-k}(y),\; x,y\in I,\; i\geq 0,\\
\gamma_{i}(ax) &=a^{i}\gamma_{i}(x),\; a\in A, \; i\geq 0,\\
\gamma_{i}(x)\gamma_{j}(x) &=\frac{(i+j)!}{i!j!}\gamma_{i+j}(x), \; x\in I,\; i,j \geq 0,\\
\gamma_{i}(\gamma_{j}(x)) &=\frac{(ij)!}{i!(j!)^{i}}\gamma_{ij}(x),\; x\in I,\; i\geq 0,\; j\geq 1.
\end{align}
The maps $\gamma_{i}$ are called the operations of divided powers. We denote a ring with divided powers by $(A,I,\gamma)$. A homomorphism of rings with divided powers  $$f: (A,I,\gamma) \rightarrow  (A',J,\delta)$$ is a ring homorphism $f: A\rightarrow B$ such that $f(I)\subset J$ and such that $f(\gamma_{i}(x))=\delta_{i}(f(x)),$ for all $i\geq 0$ and $x\in I$. A commutative algebra equipped with such a structure is called a divided powers algebra. 
\end{definition}
We remark that by relations  $(1.1)$ and $(1.4)$ we have $x^{n}=n!\gamma_{n}(x)$, where $n\in \mathbf{N}$ and $x\in I$. It follows, in prime characteristic that $x^{p}=0$, for all $x\in I$.

\textbf{Free divided powers algebras}. Let $M$ be a module over $R$, then N.Roby  constructs by Theorem $2$ in \cite{Ro2} the free divided powers algebra $\Gamma (M)$.

\begin{example}
If $k$ is a field of characteristic $0$, then any augmented commutative algebra $A$  over $k$ is an algebra with divided powers with operations $$\gamma_{i}(x):=\dfrac{x^{i}}{i!},$$ for all $x\in A_{+}$ and $i\geq 0$.
\end{example}

\begin{example}
Let $A$ be an augmented algebra over $k$ with augmentation ideal $I$. We suppose that there is an integer $m$ such that $(m-1)!$ is invertible and $I^{m}=0$. Then $A$ can be endowed with the structure of a divided powers algebra. We can define divided powers operations as follows $\gamma_{n}(x)=\dfrac{x^{n}}{n!},$  for $n\leq m-1$ and $\gamma_{n}(x)=0$, for all $n\geq m,\; x\in I$. In particular, if $k$ is a field with $char\,k=p$ and $I^{p}=0$, then $A$ becomes a divided powers algebra.
\end{example}

\begin{example}
Let $S$ be an algebra. Then $S$ is equipped with a system of divided powers $\gamma_{i}$ on the ideal $\{0\}$, where $\gamma_{0}(0)=1$ and $\gamma_{n}(0)=0$, for all $n\geq 1$. The algebra $S$ with this structure is called a trivial divided powers algebra.
\end{example}

\textbf{Special derivations}. Let $A$ be a divided powers algebra with operations $\gamma_{i}$ on an ideal $I$. A derivation $D$ of $A$ is called special if $D(\gamma_{i}(x))=\gamma_{i-1}(x)D(x)$, for all $x\in I$ and $i\geq 0$.

\subsection{Divided powers in prime characteristic} 

Let $k$ be a field of prime characteristic and let $A$ be a $k$-algebra. Following Robby in \cite{Ro2}, we call $A$ a pregraded algebra if $A$ as a $k$-module is decomposed $A=A_{0}\oplus A_{+}$, where $A_{0}$ is a subalgebra of $A$ containing $1$ and $A_{+}$ is an ideal of $A$. If $A$ is equipped with a system of divided powers operations $\gamma_{n}: A_{+}\rightarrow A_{+}$, then by a theorem of H. Cartan \cite{Car} it follows that the operations $\gamma_{n}$ are induced by the $p$-th operation $\gamma_{p}$. Moreover, let $A=A_{0}\oplus A_{+}$ be a pregraded algebra such that $x^{p}=0$ for all $x\in A_{+}$ together with a map $\pi: A_{+}\rightarrow A_{+}$ such that $\pi$ verifies the relations of the operation $\gamma_{p}$. Then it is proved by J.P. Soublin in \cite{Sou} that $A$ is equipped with a unique system of divided operators such that $\gamma_{p}=\pi$.

It follows that the category of pregraded $k$-algebras with divided powers on the ideal $A_{+}$ is equivalent to the category whose objects are pairs $(A,\pi)$, where $A$ is a pregraded algebra over $k$ and $\pi: A_{+}\rightarrow A_{+}$ is a map which we call $p$-map such that
\begin{align}
x^{p} &= 0,\;x\in A_{+}, \\
\pi(x+y) &=\pi(x)+\pi(y)+\sum_{k=1}^{k=p-1}\dfrac{(-1)^k}{k}x^{k}y^{p-k},\;x,y\in A_{+},\\
\pi(xy)  &=0,\; x,y\in A_{+},\\
\pi(a_{0} x) &=a_{0}^{p}\pi(x),\;x\in A_{+}, a_{0} \in A_{0}.
\end{align}
A homomorphism $f: (A,\pi)\rightarrow (A',\pi')$ is a ring homomorphism  $f: A\rightarrow A'$ such that $f(A_{+})\subset A'_{+}$ and such that $f\pi=\pi' f$. We denote by $pdCom$ this category and by $pdCom_{k}$  its subcategory of augmented algebras $A=k\oplus A_{+}$ with divided powers on the augmentation ideal $A_{+}$. 

Let $A\in pdCom_{k}$, then an ideal $I\subset A_{+}$  is called $p$-ideal if $\pi(I)\subset I$. The kernel $I$ of a homomorphism $f: A\rightarrow A'$ in $pdCom_{k}$ is a $p$-ideal of $A$. If $I$ is a $p$-ideal of $A$, then $A/I$ is a divided powers algebra in $pdCom_{k}$.

\begin{remark} Let $A\in pdCom$ be a divided powers algebra and $\pi,\pi'$ two $p$-maps. The map $\phi: A_{+}\rightarrow A_{+}$ defined by $\phi(x)=\pi-\pi'$ 
verifies $\phi(A_{+}^{2})=0$,
\begin{equation}
\phi(x+y)=\phi(x)+\phi(y),
\end{equation}
and
 \begin{equation}
 \phi(a_{0} x)=a_{0}^{p}\phi(x),
 \end{equation}
for all $x,y\in A_{+}$ and $a_{0} \in A_{0}$. Conversely, suppose that $(A,\pi)$is an object in $pdCom$ and $\phi: A_{+}\rightarrow A_{+}$ is a map which verifies $(1.10)$ and $(1.11)$ and such that $\phi(A_{+}^{2})=0$, then $(A,\pi')$ is a object in $pdCom$, where $\pi'=\pi+\phi$.
\end{remark}

For the classification of simple modular Lie algebras are proved theorems describing embeddings of Lie algebras into special derivations algebras. Following S. M. Skryabin and H. Strade we recall how the structure of divided powers algebra appears in Lie theory in prime characteristic. For details we refer to \cite{Skry,Stra,Ra,Ma}.

\subsection{Dual of enveloping algebra}
Let $O(n)$ be the commutative algebra defined by generators $x_{i}^{(r)}$, $i=1,...,n$, $r\geq 0$ and relations 
\begin{equation}
x_{i}^{(0)}=1,\; x_{i}^{1}=x_{i},\; x_{i}^{(s)}x_{i}^{(t)}=\dfrac{(t+s)!}{t!s!}x_{i}^{(s+t)}.
\end{equation}
We identify $x_{i}^{(1)}=x_{i}$. Let $a=(a_{1},...,a_{n})$ be a sequence, where $a_{i}$ are non negative integers. We denote $$x^{(a)}:=x_{1}^{(a_{1})}x_{2}^{(a_{2})}\cdots x_{n}^{(a_{n})}.$$
Then the elements $x^{(a)}$ with $a=(a_{1},...a_{n})$ form a basis for $O(n)$. For $j\geq 1$ we denote by $m^{j}$ the linear span of $x^{(a)}$ for all $a=(a_{1},...,a_{n})$ such that $a_{1}+...+a_{n}\geq j$. There is a decreasing chain of ideals $m^{j}$ such that $m^{i}m^{j}\subseteq m^{i+j}$. The algebra $O(n)$ is augmented $O(n)=k\oplus m^{1}$ and the ideal $m^{1}$ is maximal. Let $\hat{O}(n)$ be the completion of $O(n)$ for the corresponding topology. We denote by $d_{i}$ the derivation of $O(n)$ defined by $$d_{i}(x_{j}^{(r)})=\delta_{ij}x_{i}^{(r-1)},$$ for all $i,j=1,...,n$. We extent $d_{i}$ by continuity to a derivation of $\hat{O}(n)$. S. M. Skryabin proves by Proposition $1.6$ in \cite{Skry} that $O(n)$ (resp. $\hat{O}_{n}$) is a divided powers algebra with operations on $m^{1}$ (resp. $\hat{m}^{1}$) given by $\gamma_{r}(x_{i})=x_{i}^{(r)}$ for all $i=1,..,n$ and $n\geq 0$. The Lie algebra  $W_{n}$ of all special derivations of $O(n)$ is called Witt algebra. It is proved by Proposition $1.6$ in \cite{Skry}  that the Witt algebra $W_{n}$  (resp. $\hat{W}_{n}$) of all (continuous) special derivations of  $O_{n}$ (resp. $\hat{O}_{n}$) is free $O(n)$-module (resp. $\hat{O}(n)$-module) with basis $d_{1},...,d_{n}$.

 Let $L$ be a Lie algebra over $k$. We denote by $\hat{L}$ the $p$-envelop of $L$ and we identify $U(L)$ with $u(\hat{L})$ (for details see \cite{Stra}). We recall that the restricted enveloping $u(\hat{L})$ of $\hat{L}$ is a bialgebra with comultiplication $\Delta : u(\hat{L})\rightarrow u(\hat{L})\otimes u(\hat{L})$ given by $\Delta(x)=1\otimes x+x\otimes 1$ for all $x\in \hat{L}$ and counit $\epsilon: u(\hat{L})\rightarrow k$ with $\epsilon (x)=0$, for $x\in \hat{L}$. If $h$ is a restricted sub-algebra of $\hat{L}$,  then the set
$$F(L,h):=Hom_{u(h)}(u(\hat{L}),k)=Hom_{u(h)}(U(L),k)$$ is a commutative algebra. Indeed, using Sweedler's notation $\Delta(u)=\sum_{(u)}u_{1}\otimes u_{2}$  the  multiplication is given by $fg(u)=\sum_{(u)}f(u_{1})g(u_{2})$, where $f,g \in Hom_{u(h)}(u(\hat{L}),k)$. If $h=0$, we deduce a commutative algebra structure on the dual space $$F(L,0)=U^{*}(L).$$ 
Let $L_{0}$ be a Lie subalgebra of $L$ of codimension $n$. By Theorem $3.5$ S. M. Skryabin proves in \cite{Skry} that $F(L,\hat{L}_{0})$ is a topological divided powers algebra isomorphic to $\hat{O}(n)$.

\section{K\"ahler differentials and Quillen-Barr-Beck cohomology of divided powers algebras}

In the next subsection we recall how Quillen-Barr-Beck cohomology is defined for the category $pdCom_{k}$. For details we refer to \cite{Do}. Moreover we prove that any special derivation is a Beck derivation.
\subsection{Quillen-Barr-Beck cohomology} 

Let $A\in pdCom_{k}$ be a divided powers algebra, we denote by $(pdCom_{k}/A)_{ab}$ the category of abelian group objects of the slice category $pdCom_{k}/A$. The category $(pdCom_{k}/A)_{ab}$ is called the category of Beck $A$-modules. The forgetful functor 
$$I_{A}:(pdCom_{k}/A)_{ab} \rightarrow pdCom_{k}/A$$ admits a left adjoint functor 
$$Ab_{A}: pdCom_{k}/A \rightarrow (pdCom_{k}/A)_{ab}$$ called the abelianization functor. For $B\in pdCom_{k}/A$ and $M$ a Beck $A$-module, Quillen's cohomology groups are defined as follows

$$H^{n}(B,M)=H^{n}(Hom_{(pdCom_{k}/A)_{ab}}(\mathbf{L}Ab_{A}(B),M),$$
the right hand side of the formula denotes cohomology of a cosimplicial abelian group. 

Next, we describe Beck-derivations, Beck-modules and we determine the abelianization functor $Ab_{A}$. Let $A \in pdCom_{k}$ be a divided powers algebra and let $M$ be an $A$-module. Let $\pi: M\rightarrow M$ be a map such that
\begin{align}
\pi(m+m') &=\pi(m)+\pi(m'),\;m,m'\in M,\\
\pi(\lambda m) &=\lambda^{p}\pi(m),\;m\in M,\; \lambda\in k,\\
\pi(am) &=0,\; a\in A_{+},\;m\in M.
\end{align}
In other words, $\pi$ is a $F$-semilinear, where $F: k\rightarrow k$ is the Frobenius endomorphism and such that $\pi(am)=0,$ for all $a\in A_{+}$ and $m\in M$.

\subsection{Beck derivations}
Let $A\in pdCom_{k}$ be a divided powers algebra and let $(M,\pi)$ be a pair where $M$ is a left $A$-module, and $\pi: M\rightarrow M$ is a map which verifies $(2.1),(2.2),(2.3)$. The direct sum of $A$-modules $A\oplus M$ is ring with multiplication given by
$$(a,m)(a',m'):=(aa', am'+a'm).$$
 Then $A\oplus M$ is an augmented algebra over $k$ with augmentation ideal $$(A\oplus M)_{+}=\{(a,m),\;|\;a\in A_{+},\, m\in M \}.$$
The $p$-map $$ \pi(a,m)= (\pi(a), \pi(m)-a^{p-1}m),$$ where $a\in A_{+}$ and $m\in M$ endows $A\oplus M$ with the structure of divided powers algebra on $(A\oplus M)_{+}$. Indeed, in Lemma $3.1$ in \cite{Do} is proved that the relations $(1.6), (1.7), (1.8)$, $(1.9)$ hold. We denote by $A\oplus_{p}M$ this divided powers algebra and we obtain an object  $$A\oplus_{p}M\xrightarrow{pr_{A}} A,$$ in the comma category $pdCom_{k}/A$. For any $A\in pdCom_{k}$ we denote by $$Der_{p}(A,(M,\pi))=\{D\in Der(A,M)|\; D(\pi(a))=\pi(D(a))-a^{p-1}D(a),\,a\in A_{+}\}$$ 
This set is a group called the group of Beck derivations of $A$ by $(M,\pi)$.

\subsection{Abelian group objects} There is an isomorphism
$$Hom_{pdCom_{k}/A}(B,A\oplus_{p}M)\simeq Der_{p}(B,(M,\pi))$$ given by $$\phi \mapsto pr_{M}.$$
Therefore $A\otimes_{p}M\xrightarrow{pr_{A}} A$ is an abelian group object in the slice category $pdCom_{k}/A$. 

\begin{theorem}
Let $A\in pdCom_{k}$ be a divided powers algebra and let $B \in (pdCom_{k}/A)_{ab}$ be an abelian group object. Then there exists an isomorphism of divided powers algebras $$B\simeq A\oplus_{p} M,$$
for a specific pair $(M,\pi)$, where $M$ is an $A$-module and $\pi: M\rightarrow M$ a map which verifies the relations $(2.1),(2.2),(2.3)$.
\end{theorem}

\begin{proof}
Let $f: B\rightarrow A$  be an abelian group object in $pdCom/A$ and $M:=ker\,f$. From the group structure there is a zero morphism $z: A\rightarrow B$ which splits $f$.  Then one proves that $\omega: A\oplus_{p} M \rightarrow B$ given by $(a,m)\mapsto z(a)+m$ is a divided powers algebra isomorphism, (for details see Theorem $3.2$ in \cite{Do}).
\end{proof}

It follows that an abelian group object in $pdCom_{k}/A$ is a pair $(M,\pi)$, where $M$ is an $A$-module, and $\pi: M\rightarrow M$ is a map which verifies the relations $(2.1),(2.2),(2.3)$. Let  $k[P]$ be the polynomial ring $$k[P]=\{\sum_{i=0}^{i=n}\lambda_{i}P,\;|P\lambda_{i}=\lambda_{i}^{p}P\;\}.$$

Let $A=A_{0}\oplus A_{+}$ be a pregraded algebra with divided powers on $A_{+}$. We define the ring $V(A)$ as the ring whose underlying $k$-module is the tensor product $V(A):=A\otimes k[P]$ and the multiplication is given by:
\begin{align*}
(a\otimes 1)(a'\otimes 1)&=(aa'\otimes 1),\;\;a,a'\in A,\\
(1\otimes q)(1\otimes q')&=(1\otimes qq'),\;\; q,q'\in k[P],\\
(a\otimes 1)(1\otimes P) &=(a\otimes P),\; a\in A,\\
(1\otimes P)(a\otimes 1)&=0,\; a\in A_{+},\\
(1\otimes P)(a_{0} \otimes 1)&=(a_{0}^{p}\otimes P),\;\lambda \in A_{0}.\\
\end{align*}

\begin{remark}
Let $A,B\in pdCom$ divided powers algebras. There are embeddings of rings $i_{1}: A\rightarrow V(A)$ and $i_{2}: k[P]\rightarrow V(A)$. If $A\rightarrow B$ is a homomorphism of divided powers algebras, then a ring homomorphism $V(A)\rightarrow V(B)$ is induced. 
\end{remark}

Let $A\in pdCom_{k}$ and let $(M,\pi)$ be a Beck $A$-module. We endow $M$ with the structure of $V(A)$-module via the actions $(1\otimes P)m:=\pi(m)$, and $(a\otimes 1)m:=am$, for all $a\in A$ and $m\in M$. Conversely, any $V(A)$-module $M$ is associated to a Beck $A$-module $(\bar{M},\pi)$, where $\bar{M}$ is $M$ viewed as $A$-module via $A\rightarrow V(A)$ and  $\pi: \bar{M}\rightarrow \bar{M}$ is the map given by $$\pi(m):=(1\otimes P)m$$ for all $m\in \bar{M}$.
We easily see that the category of Beck $A$-modules is equivalent to the category $V(A)$-modules. Hereafter, we identify a $V(A)$-module with a Beck $A$-module. Under this identification, if $A\in pdCom$ is a divided powers algebra and $M$ a $V(A)$-module, we denote $Der_{p}(A,M):=Der_{p}(A, (\bar{M},\pi))$.

\begin{proposition}
Let $A\in pdCom_{k}$ be a divided powers algebra. Then any special derivation is a Beck-derivation. 
\end{proposition}

\begin{proof}
We see $A$ as an $A$-module via the multiplication. The  $A$-module $A$ together with the zero map consist  a Beck $A$-module. Let $D$ be a special derivation. Then Skryabin proved by Proposition $1.5$ (see \cite{Skry}) that $D\in Der_{k}(A)$ is special if and only if 
$$D(\gamma_{p}(x))=\gamma_{p-1}(x)D(x).$$
Since $char\,k=p$ we have $(p-1)!=(-1)\,mod(p)$. Hence,
 $$\gamma_{p-1}(x)=\frac{x^{p-1}}{(p-1)!}=-x^{p-1}.$$
It follows that $D\in Der_{p}(A,A)$.
\end{proof}

\begin{example}
We consider the augmented over $k$ algebra $A=k[X]/<X^{p}>$ with augmentation ideal $A_{+}=xA,$ where $x$ is the equivalence class of $X$ in $A$. We have $A_{+}^{p}=0$. Let  $\pi: A_{+}\rightarrow A_{+}$ be a $F$-semilinear and such that $\pi(A_{+}^{2})=0$. Then $(A,\pi)$ is a divided powers algebra, $(A_{+},\pi)$ is a Beck $A$-module, and $d\in Der_{p}(A,A_{+})$ if and only if $d(\pi(a))=\pi(d(a))$, for all $a\in A_{+}$. For example, let $k$ be a field of characteristic $3$ and  let $A=k[X]/<X^{3}>$. If $\pi: A_{+}\rightarrow A_{+}$  is the $F$-semilinear map given by $\pi(x)=x^{2}$, $\pi(x^{2})=0$, then $d(x)=2x\frac{\partial }{\partial x}\in Der_{p}(A,A_{+})$.
\end{example}

 Let $A\in pdCom_{k}$. Let  $M$ be a $V(A)$-module, we set $Pm:=(1\otimes P)m$ and $am:=(a\otimes m)$ for all $m\in M$ and $a\in A$.

\subsection{The module $\Omega_{p}(A)$ }
Let $A\in pdCom_{k}$. We define $\Omega_{p}(A)$ the $V(A)$-module with the following presentation: the generators are de symbols $da$ for $a\in A$ and the relations are
 \begin{align*}
 d(\lambda)=&0,\\
 d(\lambda a+\mu b) &=\lambda da +\mu db,\\   
               d(ab)& =adb+bda,\\
               d(\pi(c))&=Pd(c)-c^{p-1}dc,          
  \end{align*}
 where $a,b\in A$,\,$c\in A_{+}$ and $\lambda, \mu\in k$.

 \begin{theorem}
 Let $A\in pdCom_{k}$ be a divided powers algebra and let $M$ be a $V(A)$-module. Then we have the following bijection
 $$ Hom_{V(A)}(\Omega_{p}(A),M) \simeq Der_{p}(A,M).$$
 \end{theorem}

\begin{proof}
 Let $f: \Omega_{p}(A)\rightarrow M$ be a $V(A)$-homomorphism. We denote by $D_{f}: A\rightarrow \bar{M}$ the Beck derivation given by   $D_{f}(a):=f(da)$, for all $a\in A$. Then the map $f \mapsto D_{f}$ is a bijection.
 \end{proof}

Let $A\in pdCom_{k}$ a divided powers algebra. The abelianization functor $Ab_{A}: pdCom_{k}/A\rightarrow (pdCom_{k}/A)_{ab}$ is given by $$B\mapsto V(A)\otimes_{V(B)}\Omega_{p}(A).$$
 
\begin{remark}
Let $F_{p}: Sets \rightarrow pdCom_{k}$ be the free functor left adjoint to the forgetful functor $U: pdCom_{k}\rightarrow Sets$. The adjuction $(F_{p},U)$ gives rise to a cotriple $\mathbb{G}$. The cotriple resolution  associated to $\mathbb{G}$ provides a cofibrant resolution $\mathbb{G}_{*}(A)\rightarrow A$ of $A$. Therefore the cohomology groups $H^{n}(A,M)$ are the cotriple cohomology groups $H^{n}_{\mathbb{G}}(A, Der_{p}(-,M))$  of $A$ with coefficients the Beck derivations functor $Der_{p}(-,M)$. 
\end{remark} 
 
 \subsection{K\"ahler differentials }
In this subsection we prove that if $A$ is a divided powers algebra in $pdCom_{k}$, then the module of K\"ahler differentials of $A$ represents Beck derivations. First we prove the following Lemma.

\begin{lemma}
Let $A\in pdCom_{k}$ and let $\mu: A\otimes A\rightarrow A $ be the multiplication map with $I=ker\;\mu$. Then we have $$ x^{p-n}\otimes x^{n}\equiv (nx^{p-1}\otimes x)\;\;\textrm{mod}\;I^{2},\;\forall x\in I,\; 1\leq n\leq p.$$
\end{lemma}

\begin{proof}
For $n=1$ and $n=p$ is obvious. Let $0< n <p$, we will proceed by induction on $n$. We prove the claim for $n=2$. Since $(1\otimes x-x\otimes 1)^{2}\in I^{2}$ we have 
 $x^{p-2}(1\otimes x-x\otimes 1)^{2}\equiv 0\;\; \textrm{mod}\;I^{2}$. Besides,
\begin{align*}
x^{p-2}(1\otimes x-x\otimes 1)^{2} &=x^{p-2}\otimes x^{2}+x^{p}\otimes 1-2x^{p-1}\otimes x\\
                                     &=x^{p-2}\otimes x^{2}-2x^{p-1}\otimes x.
 \end{align*}
Therefore $x^{p-2}\otimes x^{2}\equiv 2x^{p-1}\otimes x\;\;\textrm{mod}\;I^{2}$. Suppose that the claim is true for all $k\leq n$. We have $$\sum_{k=0}^{k=n+1}\binom{n+1}{k}(-1)^{(n+1)-k}\, x^{(n+1)-k}\otimes x^{k}\in I^{n+1}.$$  
Therefore
\begin{align*}
0  &\equiv x^{p-(n+1)}\big(\sum_{k=0}^{k=n+1}\binom{n+1}{k}(-1)^{(n+1)-k} x^{(n+1)-k}\otimes x^{k}\big)\;\; \textrm{mod}\;I^{2}\\
   &\equiv \sum_{k=0}^{k=n+1}\binom{n+1}{k}(-1)^{(n+1)-k}\, x^{(p-k)}\otimes x^{k}\;\; \textrm{mod}\;I^{2}\\
   &\equiv \binom{n+1}{0}(-1)^{n+1} \,x^{p}\otimes 1+\sum_{k=1}^{k=n}\binom{n+1}{k}(-1)^{(n+1)-k}\,x^{p-k}\otimes x^{k}+(-1)^{0}x^{p-(n+1)}\otimes x^{n+1}\;\; \textrm{mod}\;I^{2}.\\
\end{align*}
Using the induction hypothesis we have
\begin{align*}
x^{p-(n+1)}\otimes x^{n+1} &\equiv -\sum_{k=1}^{k=n}\binom{n+1}{k}(-1)^{(n+1)-k}\,x^{p-k}\otimes x^{k}\;\; \textrm{mod}\;I^{2}\\
                           &\equiv -\sum_{k=1}^{k=n}\binom{n+1}{k}(-1)^{(n+1)-k}(k\,x^{p-1}\otimes x)\;\; \textrm{mod}\;I^{2}.\\
\end{align*}
Besides, 
$$\sum_{k=1}^{k=n} \binom{n+1}{k} (-1)^{(n+1)-k}k=-(n+1),$$
and it follows that
$$ x^{p-(n+1)}\otimes x^{n+1}=(n+1)x^{p-1}\otimes x.$$
\end{proof}

Let $A\in pdCom_{k}$ be a divided powers algebra. Then the tensor algebra $T=A\otimes A$ is an augmented algebra $T=k\oplus T_{+}$ with augmentation ideal  $$T_{+}=(k\otimes A_{+})\oplus (A_{+}\otimes k)\oplus (A_{+}\otimes A_{+}).$$
It follows from Section $8$ in \cite{Ro2} that $A\otimes A$ is a divided powers algebra in $pdCom_{k}$ with $p$-map such that $\pi(a\otimes 1)=\pi(a)\otimes 1$ and $\pi(1\otimes a)=1\otimes \pi(a),$ for all $a\in A_{+}$. The  ideal $T_{+}$ is generated by  elements of the form  $a\otimes 1$ and $1\otimes a'$ where $a, a'\in A_{+}$. The multiplication map $\mu: A\otimes A\rightarrow A $ is a homomorphism of augmented algebras and  we have $$\mu (\pi(a\otimes 1))=\pi (\mu(a\otimes 1)),$$
 $$\mu (\pi(1\otimes a'))=\pi (\mu(1\otimes a')),$$ for all $a,a'\in A_{+}$. It follows from Proposition $3$ in \cite{Ro2} that the multiplication map $\mu: A\otimes A\rightarrow A$ is a homomorphism in $pdCom_{k}$. The kernel $I$ of $\mu$ is a $p$-ideal of $A\otimes A$. The module of K\"ahler differentials of $A$ over $k$ is defined as the $A$-module $I/I^{2}$. By relation $(1.7)$ the $p$-map on $(A\otimes A)_{+}$ induce a map $\pi: I/I^{2}\rightarrow I/I^{2}$ such that $(I/I^{2},\pi)$ is a Beck $A$-module.  Moreover, we have the following proposition.

\begin{proposition}
If $A\in pdCom_{k}$ is a divided powers algebra, then the $k$-linear map $d: A \rightarrow I/I^{2}$ given by $$d(x)=(1\otimes x-x\otimes 1)+I^{2}$$ is a Beck derivation.
\end{proposition}

\begin{proof}
It is known that $d$ is a derivation.  We prove that $d$ is a Beck-derivation. For $x\in A_{+}$ we have
\begin{align*}
\pi(d(x))-\underbrace{(x\dots (x}_{p-1} d(x)\cdots) &=\\
          &=1\otimes \pi(x)-\pi(x)\otimes 1-\sum_{k=1}^{p-1}\frac{1}{k}(x^{p-k}\otimes x^{k})-x^{p-1}\otimes x+x^{p}\otimes 1.
\end{align*}
By Lemma $2.7$ we obtain that
\begin{align*}
\sum_{k=1}^{p-1}\frac{1}{k}(x^{p-k}\otimes x^{k})+x^{p-1}\otimes x &=\sum_{k=1}^{p-1}\frac{1}{k}(kx^{p-1}\otimes x)+x^{p-1}\otimes x\\
                                                                   &=(p-1)x^{p-1}\otimes x+x^{p-1}\otimes x\\
																																	 &=p(x^{p-1}\otimes x)\\
																																	 &=0.
\end{align*}
Hence, we obtain that $$d(\pi(x))=\pi(d(x))-\underbrace{(x\cdots (x}_{p-1} d(x)\cdots),$$ and $d$ is a Beck derivation.
\end{proof}
																		
\begin{theorem}
If $A\in pdCom_{k}$ and $M$ is a $V(A)$-module, then $$Hom_{V(A)}(I/I^{2},M)\simeq Der_{p}(A,M),$$ and there is an isomorphism  $$I/I^{2}\simeq \Omega_{p}(A)$$ of $V(A)$-modules.
\end{theorem}

\begin{proof}
If $f: I/I^{2}\rightarrow M$ is a $V(A)$-homomorphism, then  $f \circ d$ is a Beck derivation in  $Der_{p}(A,M)$. Indeed, by Proposition $(2.8)$ we have that $d\in Der_{p}(A,I/I^{2})$ and
\begin{align*}
(f\circ d)(\pi(x)) &=f(d(\pi(x)))\\
                   &=f(\pi(d(x)))-\underbrace{(x\cdots (x}_{p-1}d(x)\cdots)\\
									 &=\pi(f(d(x)))-\underbrace{(x\cdots (x}_{p-1}(f(d(x))\cdots).
\end{align*}
Thus we associate to each $V(A)$-homomorphism  $f: I/I^{2}\rightarrow M$, a Beck derivation in $Der_{p}(A,M)$. Conversely,  let $D\in Der_{p}(A,M)$ be a Beck derivation. We  will prove that there exists a unique $V(A)$-homomorphism  $f: I/I^{2}\rightarrow M$ such that $f\circ d=D$. We note that $I/I^{2}$ is generated by $\{d(x) | x\in A\}$ as a $V(A)$-module. If there exists such homomorphism it will be unique. Let $\phi: A\otimes A \rightarrow A\oplus_{p} M$ be the map given by $$\phi(a\otimes a')=(aa', aD(a')),$$ where $a,a'\in A$. We prove that $\phi$ is divided powers homomorphism. 
\begin{align*}
\phi((a\otimes a')(b\otimes b'))&=\phi(ab\otimes a'b')\\
                             &=(aba'b', abD(a'b'))\\
														 &=(aba'b', aba'D(b')+abb'D(a'))\\
														 &=(aa',aD(a')(bb',bD(b'))\\
														&=\phi(a\otimes a')\phi(a\otimes b').
\end{align*}
If $a,a'\in A_{+}$, then 
\begin{align*}
\pi(a\otimes a')&=\pi((a\otimes 1)(1\otimes a'))\\
                &=0.
\end{align*}
Therefore \begin{align*}
\pi(\phi(a\otimes a'))&=\pi(aa', aD(a'))\\
                     &=(\pi(aa'),\pi(aD(a'))-(aa')^{p-1}aD(a'))\\
										 &=(0,0)\\
										 &=\phi(\pi(a\otimes a')).
\end{align*}
Besides, \begin{align*}
\phi(\pi(a\otimes 1)) &=\phi(\pi(a)\otimes 1)\\
                      &=(\pi(a),\pi(a)D(1))\\
											&=(\pi(a),0),
\end{align*}
and \begin{align*}
\pi(\phi(a\otimes 1)) &=\pi(a,aD(1))\\
                      &=\pi(a,0)\\
											&=(\pi(a),0).
\end{align*}
Also, \begin{align*}
\phi(\pi(1\otimes a))          &=(\pi(a),D(\pi(a)))\\
											&=(\pi(a), PD(a)-\underbrace{(a\cdots (a}_{p-1}D(a)\cdots)\\
											&=\pi(a,D(a))\\
											&=\pi(\phi(1\otimes a)).
\end{align*}
It follows that $\phi$ is a homomorphism of divided powers algebras. Moreover, $I^{2}$ is a $p$-ideal of $A\otimes A$ and  $$\phi(I^{2})=\phi(I)\phi(I)\subset M^{2}=0$$ 
Hence, a divided powers homomorphism is induced
$$\Phi: (A\otimes A)/I^{2}\rightarrow A\oplus_{p} M.$$ 
We note that $\Phi(d(a))=(0,D(a))$. We will prove that the restriction $f=\Phi_{|_{I/I^{2}}}$ is a $V(A)$-module homomorphism $f: I/I^{2}\rightarrow M$.
We have \begin{align*}
f(\pi(d(a)) &=f\big(\pi(1\otimes a-a\otimes 1)+I^{2}\big)\\
            &=f\big((1\otimes \pi(a)-\pi(a) \otimes 1-\sum_{k=1}^{k=p-1}\frac{1}{k}a^{p-k}\otimes a^{k})+I^{2}\big)\\
						&=f\big((1\otimes \pi(a)-\pi(a)\otimes 1)+I^{2}\big)-f\big(\sum_{k=1}^{k=p-1}\frac{1}{k}a^{p-k}\otimes a^{k})+I^{2}\big)\\
						&=(0,D(\pi(a)))-\sum_{k=1}^{k=p-1}(a^{p}, \frac{1}{k} a^{p-k}D(a^{k}))\\
						&=(0,D(\pi(a)))-(0,\sum_{k=1}^{k=p-1}(\frac{1}{k} a^{p-k}ka^{k-1}D(a))\\
						&=(0,D(\pi(a)))-(0, \sum_{k=1}^{k=p-1}a^{p-1}D(a))\\
						&=(0,D(\pi(a)))-(0,(p-1)a^{p-1}D(a))\\
						&=(0,\underbrace{(a\cdots (a}_{p-1}(D(a)\cdots)+D(\pi(a)))\\
						&= (0,\underbrace{(a\cdots (a}_{p-1}(D(a)\cdots)+P\,D(a)-\underbrace{(a\cdots (a}_{p-1}D(a)\cdots)\\
						&=(0,P\,D(a))\\
						&=\pi(0,D(a))\\
						&=\pi(f(d(a)).\\
						\end{align*}
It follows that $f$ is a $V(A)$-homomorphism such that $f\circ d=D$. Finally, by Theorem $2.5$ there is an isomorphism of $V(A)$-modules $\Omega_{p}(A)\simeq I/I^{2}$. Therefore the module of K\"ahler differentials of $A$ represents  Beck derivations.
\end{proof}

\subsection{Witt algebras}

As an application of the previous section in the modular Lie theory we obtain the following corollary.

\begin{corollary} There are isomorphisms of sets
$$H^{0}(O(n),O(n))\simeq Hom_{V(O(n))}(\Omega_{p}(O(n)), O(n))\simeq W(n).$$
\end{corollary}
\begin{proof}
By Proposition $2.3$ the Witt algebra $W_{n}$ is the group of Beck derivations $Der_{p}(O(n),O(n))$. The result follows from Theorem $2.9$.
\end{proof}

\subsection{Second fundamental exact sequence}
\begin{proposition}
Let $A\in pdCom_{k}$ be a divided powers algebra and let $J$ be a $p$-ideal of $A$. There is  an exact sequence of $V(A/J)$-modules 
\begin{equation}
J/J^{2}\xrightarrow{\phi} \Omega_{p}(A) \otimes_{A} A/J\xrightarrow{\psi} \Omega_{p}(A/J)\rightarrow 0
\end{equation}

\end{proposition}

\begin{proof}

We know that there exist a short exact sequence of $A/J$-modules
\begin{equation}
J/J^{2}\xrightarrow{\phi} \Omega_{p}(A) \otimes_{A} A/J\xrightarrow{\psi} \Omega_{p}(A/J)\rightarrow 0,
\end{equation}
where  $\phi$ and $\psi$ are $A/J$-module homomorphisms given by $$\phi(x+J^{2})=d(x)\otimes (1+J),$$ and $$\psi\big(d(t)\otimes (a+J)\big)=(a+J) d(t+J),$$ 
for all $x\in J$ and $t, a\in A$. Since  $J\subset A_{+}$ by relation $(1.7)$, we obtain that the $p$-map of $A$ induce a $V(A/J)$-module structure on $J/J^{2}$. The $A/J$-module $\Omega_{p}(A)\otimes_{A}A/J$ is a $k[P]$-module via the action $$P(w\otimes (a+J))=a^{p}Pw\otimes (1+J),$$ where $\omega\in \Omega_{p}(A)$, and  $a \in A$. Hence $\Omega_{p}(A)\otimes_{A}A/J$ becomes an $V(A/J)$-module.
The map $\phi$ is a homomorphism of $k[P]$-modules. Indeed,   we see that
\begin{align*}
\phi\big(P(x+J^{2})\big) &=\phi(\pi(x)+J^{2})\\
                         &= d(\pi(x))\otimes (1+J)\\
												 &=  Pd(x)\otimes (1+J) -\underbrace{x\cdots x}_{p-1} dx\otimes (1+J)\\
												 &= Pd(x)\otimes (1+J)+0\\\
												 &=P\phi(x+J^{2}),\\
												\end{align*}
for all $x\in J$. Also, the map $\psi$ is a a homomorphism of $k[P]$-modules. In fact, we have
\begin{align*}
\psi\big(P (dt\otimes (a+J))\big) &=\psi\big(a^{p}Pdt\otimes (1+J)\big)\\
                                                    &=(a^{p}+J)\psi\big((d(\pi(t))+\underbrace{t\cdots t}_{p-1}dt)\otimes (1+J)\big)\\
                                                    &=(a^{p}+J)[\psi\big(d(\pi(t))\otimes (1+J)\big)+\psi\big(\underbrace{t\cdots t}_{p-1}dt\otimes (1+J)\big)                                               ] \\
                                                    &=(a^{p}+J)[d(\pi(t)+J)+\underbrace{(t+J)\cdots (t+J)}_{p-1}d(t+J)]\\
                                                    &=(a^{p}+J)Pd(t+J)\\
                                                    &=P (a+J)d(t+J)\\
                                                    &=P\psi\big(dt\otimes (a+J)\big),\\
                                                    \end{align*}
for all $a \in A$ and $t$. It follows that the sequence $(2.5)$ is an exact sequence of $V(A/J)$-modules. 
\end{proof}

\section{Sheaf of relative differentials of divided powers algebras}

The theory of ringed spaces with divided powers i.e PD ringed spaces is developed in \cite{Ber} by P. Berthelot. In Section $3.29$ in  \cite{Ber} is indicated that the localization of any commutative algebra with divided powers has a structure of an algebra with divided powers. In particular, for $A\in pdCom_{k}$ a divided powers algebra and $f\in A$, one can see that $A_{f}$ is a pregraded algebra equipped with a map $\pi: A_{f_+}\rightarrow A_{f_+}$ given by $\pi (\dfrac{x}{f^{i}})= \dfrac{\pi(x)}{f^{ip}}$ for all $x\in A_{+}$. We consider ringed spaces with divided powers in $pdCom$ and we give the following definition.

\begin{definition}
We call a ringed space $(X,\mathcal{O}_{X})$ a $p$-ringed space if the structure sheaf $\mathcal{O}_{X}$ is a sheaf of divided powers algebras in $pdCom$. A morphism $$(f,f^{\#}):(X,\mathcal{O}_{X})\rightarrow(Y,\mathcal{O}_{Y}) $$ of $p$-ringed spaces  is morphism of ringed spaces such that $$f^{\#}: \mathcal{O}_{Y}\rightarrow f_{*}\mathcal{O}_{X}$$ is a morphism of sheaves of divided powers algebras.
\end{definition}

\begin{example}
If $A\in pdCom_{k}$ is a divided powers algebra, then  $(Spec\,A,\mathcal{O}_{Spec\,A})$ is a $p$-ringed space. 
\end{example}

Let $A\in pdCom_{k}$ be a divided powers algebra and let  $M$ be  a $V(A)$-module. If $f\in A$, then the localization module $M_{f}=\{\dfrac{m}{f^{i}},\;\;m\in M\}$ is a $k[P]$-module via the action 
$$P (\dfrac{m}{f^{i}})=\dfrac{Pm}{f^{pi}}.$$
Thus  $M_{f}$ becomes a $V(A)$-module. 
\begin{proposition}
Let $A\in pdCom_{k}$ be a divided powers algebra  and $M$ be a $V(A)$-module. For all $f\in A$ there is a isomorphism of $V(A)$-modules
$$M_{f}\simeq M\otimes_{A} A_{f}.$$
\end{proposition}

\begin{proof}
One can easily see that $M\otimes_{A} A_{f}$ becomes a $k[P]$-module via the action 
$$P (m\otimes \dfrac{a}{f^{i}})= (a^{p}Pm \otimes \dfrac{1}{f^{pi}}),$$ where $a \in A$. It follows that $M\otimes_{A} A_{f}$ is endowed with the structure of $V(A)$-module. The map $\tau: M_{f}\rightarrow  M\otimes_{A} A_{f}$ given by $\tau(\dfrac{m}{f^{i}}):=m \otimes \dfrac{1}{f^{i}}$ is an isomorphism $V(A)$-modules.
\end{proof}
 
 \begin{remark}
The  above isomorphism of $V(A)$-modules $M\otimes  A_{f}\simeq M_{f}$ is compatible with the structure of $V(A_{f})$-modules.
\end{remark}
Let  $(X, \mathcal{O}_{X})$ be a $p$-ringed space. Next, we define the notion of $p$-$\mathcal{O}_{X}$-module.

\begin{definition}
Let $(X,\mathcal{O}_{X})$ be a $p$-ringed space. A $p$-$\mathcal{O}_{X}$ module is a sheaf of abelian groups $\mathcal{M}$ on $X$ such that $\mathcal{M}(U)$ is a $V\big(\mathcal{O}_{X}(U)\big)$-module for all $U\in Open(X)$ and the restriction maps $\mathcal{M}(U)\rightarrow \mathcal{M}(V)$ are compatible with the maps $V(\mathcal{O}_{X}(U))\rightarrow V(\mathcal{O}_{X}(V)),$ for all $V\subset U$ and $U,V\in Open(X)$. Let $\mathcal{M}$ and $\mathcal{N}$ be $p$-$\mathcal{O}_{X}$-modules.  A morphism $\phi: \mathcal{M}\rightarrow \mathcal{N}$ of sheaves of abelian groups is called a morphism of  $p$-$\mathcal{O}_{X}$-modules if $\mathcal{M}(U)\rightarrow \mathcal{N}(U)$ is a morphism of $V(\mathcal{O}_{X}(U))$-modules.  We denote by $Hom_{p-\mathcal{O}_{X}}(\mathcal{M},\mathcal{N})$ the set of $p$-$\mathcal{O}_{X}$-module morphisms.
\end{definition}

Let  $(X,\mathcal{O}_{X})$ be a $p$-ringed space and  $\mathcal{M}$ a $p$-$\mathcal{O}_{X}$ -module. By definition any $\mathcal{M}(U)$ is a $V(\mathcal (\mathcal{O}_{X}(U))$-module and one can see by Remark $2.2$ that $\mathcal{M}$ is a $\mathcal{O}_{X}$-module.

\begin{example}
Let $A\in pdCom_{k}$ be a divided powers algebra  and let $M$ be a $V(A)$-module. Then  $\widetilde{M}$ is $p$-$\mathcal{O}_{Spec\,A}$-module.
\end{example}

Let $A,B\in pdCom$ and let $\sigma: B\rightarrow A$ be a morphism of divided powers. If $G$ is a $V(B)$-module, then the tensor product $G\otimes_{B} A$ becomes a $V(A)$-module with actions $x(g\otimes y):=g\otimes xy$ and $P( g\otimes a):=a^{p}(Pg\otimes 1)$, for all $x,y  \in A, a\in A$, and $g\in G$.

\begin{lemma}
Let  $A, B\in pdCom$ and let  $\sigma: B\rightarrow A$ be a morphism of divided powers algebras. If $G$ is a $V(B)$-module and $F$ is a $V(A)$-module, then  we have a bijection
$$Hom_{V(B)}(G,F_{V(B)})\simeq Hom_{V(A)}(G\otimes_{B}A,F).$$
\end{lemma}

\begin{proof}
Let $\phi: G\rightarrow F_{V(B)}$ be a homomorphism of $V(B)$-modules. We consider the $A$-homomorphism $f_{\phi}: G\otimes_{B}A\rightarrow F$ given by
$f_{\phi}(g\otimes a)=a\phi(g)$, for all $a\in A$ and $g\in G$. For all $a \in A$ and $g\in G$
\begin{align*}
f_{\phi}(P(g\otimes a))&=f_{\phi}(a^{p}(Pg\otimes 1))\\
                           &=a^{p}\phi(Pg)\\
											&=P(a \phi(g))\\
											&=Pf_{\phi}(g\otimes a)
\end{align*}		
It follows that $f_{\phi}$ is a homomorphism of $V(A)$-modules. The map $\phi \mapsto f_{\phi}$ is a bijection.		
\end{proof}

\begin{remark}
Let $A\in pdCom_{k}$ be a divided powers algebra and $X=Spec\,A$ the spectrum of $A$. If $M$ is a $V(A)$-module and $\mathcal{F}$ is a $p$-$\mathcal{O}_{X}$-module, then we recall that there is a bijection 
$$Hom_{A}(M,\mathcal{F}(X))\simeq Hom_{\mathcal{O}_{X}}(\widetilde{M},\mathcal{F}).$$
 The restriction of this bijection on $ Hom_{V(A)}(M,\mathcal{F}(X))$ induces a bijection  $$Hom_{V(A)}(M,\mathcal{F}(X))\simeq Hom_{p-\mathcal{O}_{X}}(\widetilde{M},\mathcal{F}).$$
\end{remark}

 Let $X,Y$ be $p$-ringed spaces, $f: X\rightarrow Y$ be a morphism  and let $\mathcal{G}$ be a sheaf of divided powers in $pdCom$ on $Y$. The inverse sheaf $f^{-1}(\mathcal{G})$ is the sheaf of divided powers algebras in $pdCom$ on $X$ which is associated to the pre-sheaf of divided power algebras in $pdCom$ given by  $$U \rightarrow \underset{f(U)\subset W}{\varinjlim \mathcal{G}(W)},$$ where $W\in Open(Y)$ and $U \in Open(X)$.

\begin{proposition}
Let $X,Y$ be $p$-ring spaces and let $$(f,f^{\#}):(X,\mathcal{O}_{X})\rightarrow(Y,\mathcal{O}_{Y}) $$ be a morphism of $p$-ringed spaces and $\mathcal{G}$ a $p$-$\mathcal{O}_{Y}$-module. Then $f^{*}(\mathcal{G})$ is an $p$-$\mathcal{O}_{X}$-module and  there is an  isomorphism of sets
$$ Hom_{p-\mathcal{O}_{X}}(f^{*}(\mathcal{G}),\mathcal{F})\simeq Hom_{p-\mathcal{O}_{Y}} (\mathcal{G},f_{*}(\mathcal{F})), $$ where $\mathcal{F}$ is a $p$-$\mathcal{O}_{X}$-module.
\end{proposition}

\begin{proof}

Since $\mathcal{G}$ is $p$-$\mathcal{O}_{Y}$-module we have that  $\mathcal{G}(W)$ is a $V(\mathcal{O}_{Y}(W))$-module for all $W\in Open(Y)$. Therefore  $f^{-1}(\mathcal{G})$ is a $p$-$f^{-1}(\mathcal{O}_{Y})$-module. Besides, $\mathcal{F}$ is $p$-$\mathcal{O}_{X}$-module and $$f^{\#}: \mathcal{O}_{Y}\rightarrow f_{*}\mathcal{O}_{X}$$ is a morphism of sheaves of divided powers algebras  thus $ f_{*}(\mathcal{F})$ is a $p$-$\mathcal{O}_{Y}$-module. 

Furthermore, the morphism $f^{\#}: \mathcal{O}_{Y}\rightarrow f_{*}\mathcal{O}_{X}$ of sheaves of divided powers algebras corresponds to a morphism of sheaves of divided powers algebras $f^{-1}(\mathcal{O}_{Y})\rightarrow \mathcal{O}_{X}$. It follows that $F$ is endowed with the structure of $p$-$f^{-1}(\mathcal{O}_{Y})$-module. Moreover, we have an isomorphism of sets $\phi$

$$\phi:  Hom_{f^{-1}(\mathcal{O}_{Y})}(f^{-1}(\mathcal{G}),\mathcal{F})\rightarrow  Hom_{\mathcal{O}_{Y}} (\mathcal{G},f_{*}(\mathcal{F})). $$
The restriction of this isomorphism on  $Hom_{p-f^{-1}(\mathcal{O}_{Y})}(f^{-1}(\mathcal{G}),\mathcal{F})$ induces an isomorphism
$$ Hom_{p-f^{-1}(\mathcal{O}_{Y})}(f^{-1}(\mathcal{G}),\mathcal{F})\simeq Hom_{p-\mathcal{O}_{Y}} (\mathcal{G},f_{*}(\mathcal{F})). $$

We have that $f^{*}(\mathcal{G})=f^{-1}(\mathcal{G})\otimes_{f^{-1}(\mathcal{O}_{Y})}\mathcal{O}_{X}$ is a $p$-$\mathcal{O}_{X}$-module and by Lemma $3.7$ we obtain
$$ Hom_{p-f^{-1}(\mathcal{O}_{Y})}(f^{-1}(\mathcal{G}),\mathcal{F})\simeq Hom_{p-\mathcal{O}_{X}}(f^{-1}(\mathcal{G})\otimes_{f^{-1}(\mathcal{O}_{Y})} \mathcal{O}_{X},\mathcal{F}), $$
and  the result follows.

\end{proof}

\begin{proposition}
 Let $A,B \in pdCom$ and let $\phi: B\rightarrow A$ be a homomorphism of divided powers algebras. We denote by $X=Spec\,A$, $Y=Spec\, B$ and let $G$ be a $B$-module.  If $f: X \rightarrow Y$ is the induced $p$-ringed morphism, then there is an isomorphism of $p$-$\mathcal{O}_{X}$-modules 
$$f^{*}(\widetilde{G})\simeq \widetilde{(G\otimes_{B}A)}.$$
\end{proposition}

\begin{proof}
By Remark $3.8$ it follows that
$$ Hom_{p-\mathcal{O}_{X}}(f^{*}(\widetilde{G}),\mathcal{F})\simeq Hom_{p-\mathcal{O}_{Y}}(\widetilde{G},f_{*}(\mathcal{F}))\simeq Hom_{V(B)}(G,\mathcal{F}(X)_{V(B)}).$$
By Lemma $3.7$
$$Hom_{V(B)}(G, \mathcal{F}(X)_{V(B)})\simeq Hom_{V(A)}(G\otimes_{B}A, \mathcal{F}(X)).$$
Again by Remark $3.8$ it follows 
$$Hom_{p-\mathcal{O}_{X}}(f^{*}(\widetilde{G}),\mathcal{F})\simeq Hom_{p-\mathcal{O}_{X}}(\widetilde{G\otimes_{B}A}, \mathcal{F}).$$
By Yonada lemma follows that 
$$f^{*}(\widetilde{G})\simeq \widetilde{(G\otimes_{B}A)}.$$

\end{proof}

\subsection{Sheaf of relative differentials}

Let $A\in pdCom_{k}$ be a divided powers algebra  and $X=Spec\,A$ be an affine scheme. We consider $k$ as a trivial divided powers algebra. If  $Y=Spec\,k$, then there is a scheme morphism $f: X\rightarrow Y$ which corresponds to the divided powers algebras homomorphism $i: k\rightarrow A$. From the universal property of the tensor product it follows that the fiber product $X\times_{Y} X$ is isomorphic to $Spec\;(A\otimes_{k} A)$. The diagonal morphism 

$$\Delta: X\rightarrow X\times_{Y} X $$ 
is induced by the multiplication map $\mu: A\otimes A\rightarrow A$ which is a morphism of divided powers algebras. The sheaf of relative differentials is given by $$\Delta^{*}(\widetilde{I/I^{2}}),$$ where $I=ker\,\mu$. Since $\Delta$ is a morphism of $p$-ringed spaces it follows that $\Delta^{*}(\widetilde{I/I^{2}})$ is $p$-$\mathcal{O}_{X}$-module.

\begin{proposition}
There is  an isomorphism of $p$-$\mathcal{O}_{X}$-modules.
$$\Delta^{*}(\widetilde{I/I^{2}})\simeq (\widetilde{I/I^{2}}).$$
\end{proposition}
\begin{proof}
By Proposition $3.10$ we have
$$\Delta^{*}\big(\widetilde{I/I^{2}}\big)\simeq \widetilde{(I/I^{2}\otimes_{A\otimes A}A)} \simeq (\widetilde{I/I^{2})}.$$

\end{proof}
Let $A\in pdCom$  be a divided powers algebra with $A=A_{0}\oplus A_{+}$. Let $M$ be a $V(A)$-module we denote $$Der_{p}(A,M)=\{D\in Der(A,M)|\; D(\pi(a))=\pi(D(a))-a^{p-1}D(a),\,a\in A_{+}\}$$ and we call an element of the set $Der_{p}(A,M)$ a Beck derivation of $A$ by $M$.

\begin{definition}
Let $(X,\mathcal{O}_{X})$ be a $p$-ringed space and $\mathcal{M}$ a $p$-$\mathcal{O}_{X}$-module. A morphism of sheaves of $k$-vector spaces $D:\mathcal{O}_{X}\rightarrow \mathcal{M}$ is called a Beck derivation of $\mathcal{O}_{X}$ into $\mathcal{M}$ if for all $U\in Open(X)$ the family of $k$-linear maps
$$D_{U}: \mathcal{O}_{X}(U)\rightarrow \mathcal{M}(U),$$ is a Beck derivation i.e $D_{U}\in Der_{p}(\mathcal{O}_{X}(U), \mathcal{M}(U))$. We denote by $Der_{p}(\mathcal{O}_{X},\mathcal{M})$ the set of Beck derivations of $\mathcal{O}_{X}$ into $\mathcal{M}$.
\end{definition}

\begin{theorem}
Let $A=k\oplus A_{+}$ be a divided powers algebra and let $I$ be the kernel of the multiplication map $\mu: A\otimes A\rightarrow A$. If $X=Spec\,A$, $Y=Spec\,k$, and and $\mathcal{M}$ is a $p$-$\mathcal{O}_{X}$-module, then there is a bijection

$$Hom_{p-\mathcal{O}_{X}}(\widetilde{I/I^{2}},\mathcal{M})\simeq Der_{p}(\mathcal{O}_{X},\mathcal{M})$$ given by $a\mapsto  a \circ d_{X/Y},$
where $d_{X/Y}: \mathcal{O}_{X}\rightarrow \widetilde{I/I^{2}}$ denotes the universal derivation. 
\end{theorem}

\begin{proof}
We have seen by Theorem $2.9$ that  $\Omega_{p}(A)\simeq I/I^{2}$. The derivation $$d_{A_{f}}: A_{f}\rightarrow {\Omega_{p}(A)}_{f}$$  given by  $$d_{A_{f}}(\dfrac{a}{f^{i}})= \dfrac{1}{f^{i}} d(a)-\dfrac{1}{f^{2i}}ad(f^{i})$$ is a Beck derivation. Indeed, for all $a\in A_{+}$ we have
\begin{align*}
    d_{A_{f}}(\pi(\dfrac{a}{f^{i}})) &=d_{A_{f}}(\dfrac{\pi(a)}{f^{pi}})\\
                                    &=\dfrac{1}{f^{pi}}d(\pi(a))-\dfrac{1}{f^{2pi}}\pi(a)df^{pi}\\
                                   & =\dfrac{1}{f^{pi}}(Pd(a)-\underbrace{a\cdots a}_{p-1}d(a))+0.
\end{align*}
Moreover,
\begin{align*}
P(d_{A_{f}}(\dfrac{a}{f^{i}}))-\underbrace{\dfrac{a}{f^{i}}\cdots \dfrac{a}{f^{i}}}_{p-1}d_{A_{f}}(\dfrac{a}{f^{i}}) &= P (\dfrac{1}{f^{i}}d(a)-\dfrac{1}{f^{2i}}ad(f^{i}))- \underbrace{\dfrac{a}{f^{i}}\cdots \dfrac{a}{f^{i}}}_{p-1}\dfrac{1}{f^{i}}d(a)+0\\
                                             &=\dfrac{1}{f^{pi}}(Pd(a)-\underbrace{a\cdots a}_{p-1}d(a)).
\end{align*}
Let $M=\mathcal{M}(X)$ then there is a bijection $$Der_{p}(A,M)\simeq Der_{p}(\mathcal{O}_{X},\mathcal{M}).$$
Moreover we have a bijection $$Hom_{p-\mathcal{O}_{X}}(\widetilde{I/I^{2}},\mathcal{M})\simeq Hom_{V(A)}(I/I^{2},M).$$
Hence the theorem follows from Theorem $2.9$.
\end{proof}

\bibliographystyle{plain}

\end{document}